\documentclass[
preprint,pre,
 aps
]{revtex4-1}
\usepackage{graphicx}
\usepackage{amsmath,amssymb,amsfonts,amsthm,color,latexsym,mathrsfs}
\usepackage{amsthm,xcolor}
\newcommand\numberthis{\addtocounter{equation}{1}\tag{\theequation}}

\newtheorem{theorem}{Theorem}[section]
\newtheorem{lemma}[theorem]{Lemma}
\newtheorem{coro}[theorem]{Corollary}

\def\half{\frac{1}{2}}
\def\veps{\varepsilon}

\def\cal{\mathcal}

\def\Dt{\Delta t}
\def\dt{\delta t}
\def\wt{\widetilde}
\def\wh{\widehat}
\def\sL{\mathscr{L}}
\def\sA{\mathscr{A}}
\def\sB{\mathscr{B}}
\def\sR{\mathscr{R}}
\def\sD{\mathscr{D}}
\def\sI{\mathscr{I}}

\newcommand{\RNum}[1]{\uppercase\expandafter{\romannumeral #1\relax}}
\def\bi{{\bf I }}
\def\bii{{\bf II }}
\def\biii{{\bf III}}

\begin{document}

\preprint{APS/123-QED}

\title{Some New Symplectic Multiple Timestepping Methods for Multiscale Molecular Dynamics Models}

\author{Chao Liang}
\email{charlesleung2009@gmail.com}
%\altaffiliation[Currently at: ]{Morgan Stanley, New York.}
\author{Xiaolan Yuan}
\email{xxy5111@psu.edu}
\altaffiliation[Currently at: ]{School of Mathematics and Statistics, Wuhan University, Wuhan, Hubei, 430072, China.}
\author{Xiantao Li}
\email{xli@math.psu.edu}
\affiliation{Department of Mathematics,\\ The Pennsylvania State University, University Park, PA 16802.}

\begin{abstract}
 We derived a number of numerical methods to treat biomolecular systems with multiple time scales. Based on the splitting of the operators associated with the slow-varying and fast-varying forces,  new multiple time-stepping (MTS) methods are obtained by eliminating the dominant terms in the error. These new methods can be viewed as a generalization of the impulse method \cite{grubmuller1991generalized,tuckerman1992reversible}. In the implementation of these methods, the long-range forces only need to be computed on the slow time scale, which reduces the computational cost considerably.  Preliminary analysis for the energy conservation property is provided. 
\end{abstract}

\date{\today}


%\begin{document}
\maketitle

\section{Introduction}

Simulations based on dynamics models for bio-molecular systems have become standard computational tools for studying small-scale motions or calibrating coarse-grained models. A core element in these simulations is the
numerical integrator for solving the underlying differential equations. A well known challenge is the presence of multiple time scales, usually arising from the various types of molecular interactions.  For example, the bond length and bond angle contributions constitute the most dominant terms, and they determine the fastest time scales. As a result, the step size of a numerical integrator has to be selected accordingly. For instance, the step size for the Verlet's method needs to satisfy \cite{schlick2010molecular}, 
\begin{equation}
 \delta t \le \frac{2}{\omega_\text{max}},
\end{equation}
where $\omega_\text{max}$ is the maximum frequency.  The bound on the step size is typically very small (femto-seconds or less).

In practice, the small step size $\delta t$ imposes a significant limitation on how long the simulations can be conducted. In particular, the most expensive part of the computation is the force calculation. While interactions via changes of bond length and bond angles are short-ranged and easy to compute, there are long-range interactions (electrostatic) that take up considerable CPU time.

A remarkable approach to overcome this difficulty is the multiple time stepping (MTS) method,
referred to as Verlet-I in  \cite{grubmuller1991generalized} and r-RESPA in  \cite{tuckerman1992reversible}. 
In sharp contrast to conventional integrators, this method involves multiple step sizes, which correspond to
the time scales that forces of different nature determine.  The MTS method is also known as the impulse method. The implementation of the impulse method is quite simple: For each small time step $\delta t$, one integrates the ODEs with only the fast force, and for each large time step \(\Dt\), one updates the momentum by applying the slow force, which is often described as a half-step `kicking', `oscillating', and another half-step `kicking'. The method can be mathematically written as a symmetric splitting of the Liouville operator  \cite{tuckerman1992reversible}, and it has been implemented in several software packages  (e.g., TINKER \cite{ponder1987tinker}). From a practical viewpoint, the impulse method eliminates the need to evaluate the
long-range forces at every step with step size $\dt$, which is the most expensive part of the dynamics simulations. As a result, the computation is sped up considerably. Another important mathematical property is the symplectic structure, which for Hamiltonian systems, is critical to preserve the total energy and produce the correct statistics.  

Despite the popularity of the impulse method, these are known practical issues. For instance, instability has been observed for some particular choices  of $\Dt$ \cite{schlick2010molecular,schlick1999some,barth1998extrapolation}. Furthermore, the order of accuracy is largely unknown. 
This has motivated several  generalization and improvement of the impulse method. For instance, it is generalized to incorporate molecular interactions of multiple scales ($\ge 3$) by Procacci and Marchi \cite{procacci1996taming}, in which the evolution operator is split into operators that represent the forces of different magnitude. The LN method \cite{barth1998overcoming} combines the slow and fast forces via extrapolation, and Langevin dynamics that represents a heat bath is introduced to reduce energy drift. Another important development is the MOLLY mothod  \cite{garcia1998long,izaguirre1999longer}, where the impulse method is modified by properly averaging the slow force.  The goal of the LN and MOLLY methods has been to overcome the resonance instability of the impulse method, which arises when the slow time step $\Dt$ is a multiple of half of the period of the fast mode.

This work continues with the operator-splitting procedure that was used to derive the original impulse method \cite{tuckerman1992reversible}. 
More specifically, we seek splitting methods that involve more fractional steps. In order to determine the coefficients of the splitting methods, a high-order error expansion is needed, which in general, is a tremendous challenge, especially when multiple fractional steps are involved. Based on the symbolic code that we had recently developed \cite{Liang2015}, we are able to obtain the error up to any prescribed order. This constitutes the basis for determining the coefficients to maximize the order of accuracy. 

The different magnitude of the fast and slow forces has motivated us to introduce another parameter $\veps$, which indicates the resulting time scale separation. Our analysis reveals that some of the terms in the error 
are weighted by large factors such as  $\frac{1}{\veps^2}$. This prompts us to re-consider their role in the error. By eliminating dominating terms in the error, we obtain three new impulse methods, which involve two to four fractional steps. Numerical tests suggest that the new methods have improved accuracy. Furthermore, these methods also have better performance in conserving the total energy. These numerical observations can be interpreted with some preliminary analysis.

%Another challenge comes from the multiscale nature of the dynamical system. In traditional numerical methods %for ODEs, one often arrange the error based on the power of the step size, e.g., $\Delta t$. The current %problem, however, contains another parameter,  $\veps$, that represents the  fast time scale. The presence of %the two small parameters has an important influence on the competition of the terms in the error. 

\section{Mathematical Formulations}\label{sec: math}

In this section, we present the mathematical framework for deriving new impulse methods. 

\subsection{Problem setup}
To explicitly indicate the fast-varying force in the dynamics, we consider a dynamical system in the following form,
\begin{equation}\label{eq:hamilton}
\left\{
\begin{aligned}
&\dot{x_i}=v_i,\\
&\dot{v_i}=\frac{1}{\veps^2}g_i(x)+f_i(x), \quad i=1,2,\cdots,N.
\end{aligned}
\right.
\end{equation}
In particular, the fast time scale is indicated by the parameter $\veps,$ with $\veps \ll 1.$ 
%In the numerical implementation, we should have $\delta t=\cal{O}(\veps).$  

Numerical methods for such ODEs, especially the one-step methods, are based on the approximation of the evolution operator, which for this autonomous system can be expressed as,
\begin{equation}\label{eq: semi}
 \big(x(t),v(t)\big)= e^{t\mathscr{L}}\big(x,v\big).
\end{equation}
Here $(x,v)$ denotes the solution at $t=0.$ Furthermore,
the differential operator $\sL$ is defined as,
\begin{equation}
  \mathscr{L}=\sum_i v_i \partial_{x_i}+\sum_i (\frac{1}{\veps^2}g_i + f_i) \partial_{v_i}. 
\end{equation}
All the partial derivatives are defined with respect to the initial data $(x,v)$.

An important class of methods have been based on the splitting of the operator $\sL$ in the semi-group operator \eqref{eq: semi}.  As a one-step method, it suffices to consider the approximation formulas for $e^{\Dt \sL},$ since the same formula can be applied to all  the following steps.
For the present problem, we  define
\begin{equation}\label{eq: L1L2}
\left\{
\begin{aligned}
\mathscr{L}_1=& \mathscr{L}_1^x+\frac{1}{\veps^2}\mathscr{L}_1^v,\\
\mathscr{L}_1^x=& \sum_i  v_i \partial_{x_i}\\
\mathscr{L}_1^v=&\frac{1}{\veps^2} \sum_i  g_i\partial_{v_i}\\ 
\mathscr{L}_2=&\sum_j  f_j \partial_{v_j}.\\
\end{aligned}\right.
\end{equation}

The popular impulse method is based on the following specific splitting scheme \cite{tuckerman1992reversible}, 
\begin{equation}\label{eq: imp0}
  e^{\Delta t\sL} \approx e^{\half\Delta  t \mathscr{L}_2}e^{\Delta  t \mathscr{L}_1} e^{\half\Delta  t\mathscr{L}_2}.
\end{equation}

In this paper, we seek splitting methods of the general form,
\begin{equation}\label{eq: mts}
e^{\Delta t\sL} \approx \prod_{i=1}^{k} e^{c_i\Delta  t\mathscr{L}_2}e^{d_i\Delta  t\mathscr{L}_1}
\end{equation}
Here the right hand side is interpreted as,
\[
  e^{c_1\Delta  t\mathscr{L}_2}e^{d_1\Delta  t\mathscr{L}_1} e^{c_2\Delta  t\mathscr{L}_2}e^{d_2\Delta  t\mathscr{L}_1}  \cdots e^{c_k\Delta  t\mathscr{L}_2}e^{d_k\Delta  t\mathscr{L}_1}.
\]
For example, in the impulse method \eqref{eq: imp0}, we have $c_1 = c_2 = \half, d_1 = 1$ and $d_2 = 0.$

Notice that the step $e^{c_i \Delta t \mathscr{L}_2}$  can be implemented exactly, since the corresponding differential equations have explicit solutions. On the other hand, the step  $e^{c_i \Delta t \mathscr{L}_1}$  requires further approximation. The corresponding differential equations are,
\begin{equation}\label{eq:hamilton1}
\left\{
\begin{aligned}
&\dot{x_i}=v_i,\\
&\dot{v_i}=\frac{1}{\veps^2}g_i(x).
\end{aligned}
\right.
\end{equation}
For instance, it can be approximated by the Verlet's method with  step size $\delta t \ll \veps$. Namely,
\begin{equation}
 e^{c_i \Dt \sL_1} \approx \Big[ e^{\frac{\dt}{2}\sL_1^v } e^{\dt \sL_1^x}  e^{\frac{\dt}{2}\sL_1^v } \Big]^m,
\end{equation}
in which $\dt=\frac{c_i \Dt}{m}$. In this case, the error is expected to be on the order of $\cal{O}(\delta t^2)$, which is negligible compared to the larger time step $\Delta t$. Therefore, we will simply assume that the step $e^{c_i \Delta t \mathscr{L}_1}$ is implemented exactly, and focus primarily on the error from the splitting at the larger steps $\Delta t.$ 

\smallskip
As discussed in the introduction, the advantage of introducing such splitting methods arises when the fast force $g$ is short-ranged, the computation of which is much less expensive than that of $f(x).$ Therefore, compared to a direct discretization of the ODEs \eqref{eq:hamilton}, where $f(x)$ is computed at {\it every} step, the splitting methods can be implemented at a reduced cost. Furthermore, one can show that all such splitting methods are symplectic, when the ODEs form a Hamiltonian system \cite{hairer2006geometric}. Finally, the computer implementation of these methods is quite straightforward, as demonstrated in \cite{FrSm02,martyna1996explicit,tuckerman1992reversible}. 
\smallskip

On the other hand, the accuracy of this method has not been extensively studied. Our mathematical formulation relies on an explicit representation of the error for the splitting method \eqref{eq: mts}. More specifically,
given the two differential operators $\mathscr{L}_1$ and $\mathscr{L}_2$, the number of fractional steps $k$, and the coefficients $(c_1,c_2,\cdots c_k)$ and $(d_1,d_2,\cdots, d_k)$, we seek $\widetilde{\mathscr{R}}$, such that:
\begin{equation}\label{productSplitting}
\prod_{i=1}^{k} e^{c_i\Delta  t \mathscr{L}_2}e^{d_i\Delta  t\mathscr{L}_1}=\widetilde{\mathscr{R}}(c_0,\cdots c_k;d_0,\cdots d_k)e^{\Delta  t\mathscr{L}}
\end{equation}
where the near-identity operator $\widetilde{\sR}$ can be expanded as follows, 
\begin{equation}\label{eq: R}
\widetilde{\mathscr{R}}(c_0,\cdots c_k;d_0,\cdots d_k)-\sI=\widetilde{\mathscr{Z}}_1{\Delta  t}+\widetilde{\mathscr{Z}}_2{\Delta  t}^2+\cdots \widetilde{\mathscr{Z}}_n{\Delta  t}^n+\cdots.
\end{equation}
Here $\sI$ is the identity operator, and $\widetilde{\mathscr{Z}}_j$ contains terms with $j$-folded operator multiplications, with coefficients represented by $(c_1, \cdots c_k; d_1, \cdots, d_k)$, i.e., $c_1c_2d_2 \mathscr{L}_2\mathscr{L}_2\mathscr{L}_2\mathscr{L}_2 \in \wt{\mathscr{Z}}_4$. All terms appeared in $\widetilde{\mathscr{Z}}_j$ are linearly independent. With the same effort, one can also express \eqref{eq: R}
in the following form,
\begin{equation}
\prod_{i=1}^{k} e^{c_i\Delta  t \mathscr{L}_2}e^{d_i\Delta  t\mathscr{L}_1}=e^{\Delta  t\mathscr{L}}\wh{\mathscr{R}}(c_0,\cdots c_k;d_0,\cdots d_k).
\end{equation}
The operator $\wh{\sR}$ can be expanded in a similar format.

The local accuracy of the splitting method is determined by the magnitude of the first few none-zero operator coefficients on the right hand side of \eqref{eq: R}.  Therefore, it is necessary to derive the explicit expressions for the first few coefficients. This procedure will be described in the next section.

\subsection{Finding explicit forms of $\wt{\sR}$}

Without loss of generality, let us consider two operators $\sA$ and $\sB$. 
The starting point of our analysis is the following expansion,
\begin{equation}
e^\sA e^\sB = {\mathscr{R}}_{(\sA,\sB)}e^{\sA+\sB}, 
\end{equation}
or alternatively,
\begin{equation}
\mathscr{R}_{(\sA,\sB)} = e^\sA e^\sB e^{-\sA-\sB}. 
\end{equation}

We want to express this approximation in operator multiplication form, instead of the exponential form in \cite{yoshida1990construction}, because multiplication form often gives back cleaner expressions \cite{Liang2015}. Direct computation shows that:
\begin{equation}\label{myExpansion}
\begin{aligned}
\mathscr{R}_{(\sA,\sB)} = & \frac{1}{2}(\sA\sB-\sB\sA) \\
 &+ \frac{1}{6}(2\sA^2\sB-4\sA\sB\sA-\sA\sB^2+2\sB\sA^2+2\sB\sA\sB-\sB^2\sA)\\
 & +\cdots
\end{aligned}
\end{equation}
More terms are available, simply by Taylor expansion of each term.

We now consider further splitting of $\sA + \sB$, in the general form of,
\begin{equation}\label{eq: splitAB}
\prod_{i=1}^{k} e^{c_i \sA}e^{d_i\sB}.
\end{equation}
 
The key to obtain the expansion of the error for \eqref{eq: splitAB} is to repeatedly use the above formula:
\begin{align*}
&\prod_{i=1}^k e^{c_i\sA}e^{d_i\sB} \\
= & \sR_{(c_1 \sA,d_1 \sB)}e^{c_1\sA+d_1\sB} \prod_{i=2}^k e^{c_i\sA}e^{d_i\sB} \\
=&  \sR_{(c_1 \sA, d_1\sB)}  \sR_{(c_1\sA+d_1\sB, c_2\sA)} e^{c_1\sA+d_1\sB+c_2\sA} e^{d_2\sB} \prod_{i=3}^k 
 e^{c_i\sA}e^{d_i\sB} \\
=&   \sR_{(c_1\sA,d_1\sB)}  \sR_{(c_1\sA+d_1\sB, c_2\sA)}  \sR_{( c_1\sA+c_2\sA+d_1\sB, d_2\sB)} e^{c_1\sA+c_2\sA+d_1\sB+d_2\sB} 
  \prod_{i=3}^k e^{c_i\sA}e^{d_i\sB}\\
=&  \cdots \\
=&  \prod_{i=1}^k \sR_{ (C_i\sA+D_{i-1}\sB, d_i\sB)}  \sR_{(C_i\sA+D_i\sB, c_{i+1}\sB)} e^{C_k \sA+ D_k \sB}\\
\triangleq & \wt{\sR}(c_1,\cdots c_k; d_1, \cdots, d_k) e^{C_k\sA+D_k \sB}
 \numberthis \label{eq:myApprox}
\end{align*}
where $$C_i = \sum_{j=1}^i c_j, \quad D_i = \sum_{j=1}^i d_j,$$ for $1 \leq i \leq k$, 
are the partial sums of the coefficients, and we define $C_0 = D_0 = 0$. 

This systematic procedure will be applied to analyze the splitting methods \eqref{eq: mts} for the multiscale ODEs \eqref{eq:hamilton}. In particular, we let $\sA=\Dt\sL_2$ and $\sB=\Dt\sL_1$. We require that  $C_k=D_k=1$, which is the consistency condition for one-step methods \cite{deuflhard2002scientific}. This leads to the error for \eqref{productSplitting} and \eqref{eq: R}. In particular, we now have,
\begin{equation}\label{eq: Rt}
\begin{aligned}
&\wt{\sR}(c_1,\cdots c_k; d_1, \cdots, d_k)\\
=&  \prod_{i=1}^k \sR_{\big( \Dt(C_i\sL_2+D_{i-1}\sL_1), \Dt d_i\sL_1\big)}  \sR_{\big(\Dt (C_i\sL_2+D_i\sL_1), \Dt c_{i+1}\sL_1\big)}.
\end{aligned} 
\end{equation}
Clearly, it would be a lengthy procedure to carry out the multiplication of these operators. Fortunately, we have developed a symbolic code \cite{Liang2015} to obtain the expansion of the error. This constitutes the basis to examine the coefficients of the error, from which the order of the accuracy can be determined, controlled and improved. In the next section, we discuss numerous cases.

\medskip

\subsection{The selection of the coefficients $c_i$ and $d_i$ based on the error expansion}

We now return to integrators with the general form \eqref{eq: mts}. We will examine the cases $k=2$, $k=3$ and $k=4$ separately.

\subsubsection{Error expansion for the case $k=2.$}

 We first consider the case $k=2.$
 Based on the analysis from the previous section,  we found,
\begin{equation}\label{eq: dt23}
 \begin{aligned}
&\wt{\mathscr{R}}(c_1,c_2=1-c_1;d_1,d_2=1-d_1) - \sI\\
=&((c_1-1)d_1+\half)(\mathscr{L}_2\mathscr{L}_1-\mathscr{L}_1\mathscr{L}_2)\Delta t^2\\
&+(\frac{d_1^2(1-c_1)}{2}-\frac{1}{6})(\mathscr{L}_1\mathscr{L}_1\mathscr{L}_2-2\mathscr{L}_1\mathscr{L}_2\mathscr{L}_1+\mathscr{L}_2\mathscr{L}_1\mathscr{L}_1)\Delta t^3\\
& +(\frac{d_1(c_1^2-1)}{2}+\frac{1}{3})(\mathscr{L}_1\mathscr{L}_2\mathscr{L}_2-2\mathscr{L}_2\mathscr{L}_1\mathscr{L}_2+\mathscr{L}_2\mathscr{L}_2\mathscr{L}_1)\Delta t^3+ \ldots
 \end{aligned}
\end{equation}

In the traditional impulse method, $c_1=c_2=\frac12$, and it may appear as if the error is $\mathcal{O}(\Delta t^3)$.  
But equation \eqref{eq: L1L2} suggests that $\mathscr{L}_1^v$ is $\mathcal{O}(\frac{1}{\veps^2})$, which clearly indicates that a further inspection of the order of the error terms is needed. To this end, let us define,
 \begin{equation}\label{eq: D23}
\begin{aligned}
 \mathscr {D}_{21}\triangleq &\mathscr{L}_2\mathscr{L}_1-\mathscr{L}_1\mathscr{L}_2,\\
 \mathscr D_{31}\triangleq& \mathscr{L}_1\mathscr{L}_1\mathscr{L}_2-2\mathscr{L}_1\mathscr{L}_2\mathscr{L}_1+\mathscr{L}_2\mathscr{L}_1\mathscr{L}_1,\\
\mathscr D_{32}\triangleq &\mathscr{L}_1\mathscr{L}_2\mathscr{L}_2-2\mathscr{L}_2\mathscr{L}_1\mathscr{L}_2+\mathscr{L}_2\mathscr{L}_2\mathscr{L}_1.
\end{aligned} 
\end{equation}

We first begin with the observation that,
\begin{lemma}\label{eq:lem1}
The following identity holds,
\begin{equation}
\mathscr{L}_2\mathscr{L}_1^v=\mathscr{L}_1^v\mathscr{L}_2.
\end{equation}
\end{lemma}
The equality can be checked directly,
\begin{equation}\label{eq:L1vL2}
\mathscr{L}_2\mathscr{L}_1^v-\mathscr{L}_1^v\mathscr{L}_2=\sum_i \sum_j f_jg_i(\partial_{v_j}\partial_{v_i}-\partial_{v_i}\partial_{v_j}).
\end{equation}
By assuming that the flow of the ODEs is sufficiently smooth so that $\partial_{v_j,v_i}^2=\partial_{v_i,v_j}^2$, we will get $\mathscr{L}_2\mathscr{L}_1^v=\mathscr{L}_1^v\mathscr{L}_2.$

Following this calculation, we find that,
\begin{theorem}\label{eq:thm0}
The order of the operators in \eqref{eq: D23} is given as follows,
\begin{equation}
\begin{aligned}
&1)\mathscr D_{21}=\cal{O}(1),\\
&2)\mathscr D_{31}=\cal{O}(\frac{1}{\veps^2}),\\
&3)\mathscr D_{32}=\cal{O}(1).
\end{aligned}
\end{equation}
\end{theorem}

We briefly outline the calculation here: \\
1)
$\mathscr{L}_2\mathscr{L}_1=\mathscr{L}_2(\mathscr{L}_1^x+\frac{1}{\veps^2}\mathscr{L}_1^v)$ and $ \mathscr{L}_1\mathscr{L}_2=(\mathscr{L}_1^x+\frac{1}{\veps^2}\mathscr{L}_1^v)\mathscr{L}_2$. So $\mathscr{L}_2\mathscr{L}_1-\mathscr{L}_1\mathscr{L}_2=\mathscr{L}_2\mathscr{L}_1^x-\mathscr{L}_1^x\mathscr{L}_2+\frac{1}{\veps^2}(\mathscr{L}_2\mathscr{L}_1^v-\mathscr{L}_1^v\mathscr{L}_2)=\mathscr{L}_2\mathscr{L}_1^x-\mathscr{L}_1^x\mathscr{L}_2$ (By Lemma \ref{eq:lem1}).\\
Further computation also shows that, 
\begin{equation}\label{eq:D2}
\mathscr D_{21}=\mathscr{L}_2\mathscr{L}_1^x-\mathscr{L}_1^x\mathscr{L}_2=
 \sum_i  f_i \partial_{x_i} - \sum_i \sum_j v_i  \partial_{x_i} f_j  \partial_{v_j} \neq 0.
\end{equation}
Hence $\mathscr{L}_2\mathscr{L}_1^x\neq \mathscr{L}_1^x\mathscr{L}_2$, i.e., $\mathscr{L}_1^x$ and $\mathscr{L}_2$ in general do not commute. Therefore,
 $\mathscr  D_{21}=\cal{O}(1)$\\
\\
2)
By direct computation we obtain that,
\begin{equation}\label{eq:t3}
\begin{aligned}
&\mathscr{L}_1\mathscr{L}_1\mathscr{L}_2-2\mathscr{L}_1\mathscr{L}_2\mathscr{L}_1+\mathscr{L}_2\mathscr{L}_1\mathscr{L}_1=\mathscr{L}_1^x\mathscr{L}_1^x\mathscr{L}_2+\mathscr{L}_2\mathscr{L}_1^x\mathscr{L}_1^x-2\mathscr{L}_1^x\mathscr{L}_2\mathscr{L}_1^x\\
&+\frac{1}{\veps^2}(\mathscr{L}_1^x\mathscr{L}_1^v\mathscr{L}_2+\mathscr{L}_1^v\mathscr{L}_1^x\mathscr{L}_2+\mathscr{L}_2\mathscr{L}_1^x\mathscr{L}_1^v+\mathscr{L}_2\mathscr{L}_1^v\mathscr{L}_1^x-2\mathscr{L}_1^x\mathscr{L}_2\mathscr{L}_1^v-2\mathscr{L}_1^v\mathscr{L}_2\mathscr{L}_1^x)\\
&+\frac{1}{\veps^4}(\mathscr{L}_1^v\mathscr{L}_1^v\mathscr{L}_2+\sL_2\mathscr{L}_1^v\mathscr{L}_1^v-2\mathscr{L}_1^v\mathscr{L}_2\mathscr{L}_1^v)
\end{aligned}
\end{equation}
From  Lemma \ref{eq:lem1}, we deduce that $\mathscr{L}_1^v\mathscr{L}_1^v\mathscr{L}_2+\sL_2\mathscr{L}_1^v\mathscr{L}_1^v-2\mathscr{L}_1^v\mathscr{L}_2\mathscr{L}_1^v=0$. 
The $\cal{O}(\frac{1}{\veps^2})$ term can be simplified to $\sD_{21}\sL_1^v-\sL_1^v\sD_{21}$, but it  is nonzero in general.

3) Following a similar calculation, we  get,
\begin{equation}\label{eq:D32}
\begin{aligned}
&\mathscr D_{32}=\mathscr{L}_1\mathscr{L}_2\mathscr{L}_2-2\mathscr{L}_2\mathscr{L}_1\mathscr{L}_2+\mathscr{L}_2\mathscr{L}_2\mathscr{L}_1\\
&=\mathscr{L}_1^x\mathscr{L}_2\mathscr{L}_2-2\mathscr{L}_2\mathscr{L}_1^x\mathscr{L}_2+\mathscr{L}_2\mathscr{L}_2\mathscr{L}_1^x+\frac{1}{\veps^2}(\mathscr{L}_1^v\mathscr{L}_2\mathscr{L}_2-2\mathscr{L}_2\mathscr{L}_1^v\mathscr{L}_2+\mathscr{L}_2\mathscr{L}_2\mathscr{L}_1^v)\\
&=\mathscr{L}_1^x\mathscr{L}_2\mathscr{L}_2-2\mathscr{L}_2\mathscr{L}_1^x\mathscr{L}_2+\mathscr{L}_2\mathscr{L}_2\mathscr{L}_1^x\\
&= \sD_{21}\sL_2-\sL_2\sD_{21} =\cal{O}(1)
\end{aligned}
\end{equation}

Based on these explicit estimates, we now have for  the original impulse method \eqref{eq: imp0}:
\begin{coro}
 When $c_1=c_2=\half, d_1=1, d_2=0$, we have 
 \begin{equation}
 \widetilde{\sR}(c_1,c_2;d_1,d_2)-\sI=\mathcal{O}(\frac{\Delta t^3}{\veps^2}) + \cdots
\end{equation}
\end{coro} 
 
\medskip 
 
 The important observation in this analysis is that the term $\sD_{31}$ contains a large factor (\(\frac{1}{\veps^2}\)). 
This motivates a different choice of the parameter: 
\begin{coro}
In the case when  
\begin{equation}\label{cor: imp1}
c_1=\frac{1}{4}, c_2=\frac{3}{4}, d_1=\frac{2}{3}, d_2=\frac{1}{3},
\end{equation}
we have,
 \begin{equation}
  \widetilde{\sR}(c_1,c_2;d_1,d_2)-\sI=\cal{O}(\Delta t^3) + \cdots.
\end{equation}
\end{coro}
In this case, we have abandoned the symmetry of the method, and chosen the parameters to eliminate
the terms of the order $\cal{O}({\Delta t^2})$ and $\mathcal{O}(\frac{\Delta t^3}{\veps^2})$ altogether. For better reference, we will call the original impulse method \eqref{eq: imp0} the {\bf impulse I} and the non-symmetric method \eqref{cor: imp1} {\bf impulse II}.

In the case when $k=2$, it is clear
 that the original impulse method is the only symmetric method. In order to
obtain other symmetric methods, we need to consider splitting methods with more fractional steps. This will be discussed further in the next two sections.

\subsubsection{ Expansions of the error for $k=3$}

 When $k=3$, we can choose $d_1=d_2=\half,$ $d_3=0,$ $c_3=c_1,$ and $c_2=1-2c_2$ to form a symmetric integrator.
 In this case,  the error  $\widetilde{\sR}$  is given by,
\begin{equation}\label{eq: k=3}
\begin{aligned}
\widetilde{\sR}&-\sI=(\frac{c_1}{4}-\frac{1}{24})(\mathscr{L}_1\mathscr{L}_1\mathscr{L}_2-2\mathscr{L}_1\mathscr{L}_2\mathscr{L}_1+\mathscr{L}_2\mathscr{L}_1\mathscr{L}_1)\Delta t^2\\
&+(\frac{c_1^2}{2}-\frac{c_1}{2}+\frac{1}{12})(\mathscr{L}_1\mathscr{L}_2\mathscr{L}_2-2\mathscr{L}_2\mathscr{L}_1\mathscr{L}_2+\mathscr{L}_2\mathscr{L}_2\mathscr{L}_1)\Delta t^3+\cdots\\
&=(\frac{c_1}{4}-\frac{1}{24})\mathscr D_{31}\Delta t^3+(\frac{c_1^2}{2}-\frac{c_1}{2}+\frac{1}{12}) \mathscr D_{32}\Delta t^3+\cdots
\end{aligned}
\end{equation}

Based on the estimate in theorem \ref{eq:thm0}, we choose to eliminate the term $\sD_{31}$, yielding,
\begin{coro}
  When $c_1=c_3=\frac{1}{6}$, $c_2=\frac{2}{3}$, $d_1=d_2=\half$, $d_3=0$, we have $\widetilde{\sR}-\sI=\frac{1}{72} \sD_{32} \Delta t^3 + \cdots.$ 
\end{coro}  

This method will be referred to as {\bf impulse} \biii. For $k=3$, there are also non-symmetric methods. But we will continue to consider the case $k=4$.

\subsubsection{ Expansions of the error for $k=4$}
Finally, we will further explore the operator-splitting methods for the case $k=4$. 
In this case, a symmetric method can be constructed by choosing $c_1,c_2=\frac{1}{2}-c_1,c_3=\frac{1}{2}-c_1,c_4=c_1;d_1,d_2=1-2d_1,d_3=d_1,d_4=0$. Up to $\Dt^4$ terms, we have,
\begin{equation}\label{eq: dt34}
 \begin{aligned}
&\widetilde{\sR}(c_1,c_2,c_3,c_4,d_1,d_2,d_3,d_4) - \sI\\
=&(c_1d_1-\half{d_1}-c_1d_1^2+\half{d_1^2}+\frac{1}{12})(\mathscr{L}_1\mathscr{L}_1\mathscr{L}_2-2\mathscr{L}_1\mathscr{L}_2\mathscr{L}_1+\mathscr{L}_2\mathscr{L}_1\mathscr{L}_1)\Delta  t^3\\
& +(\frac{1}{4}d_1-c_1d_1+c_1^2d_1-\frac{1}{24})(\mathscr{L}_1\mathscr{L}_2\mathscr{L}_2-2\mathscr{L}_2\mathscr{L}_1\mathscr{L}_2+\mathscr{L}_2\mathscr{L}_2\mathscr{L}_1)\Delta  t^3\\
&+(\half{c_1d_1}-\frac{1}{4}d_1-\half{c_1d_1^2}+\frac{1}{4}d_1^2+\frac{1}{24})\\
&\quad \times (\mathscr{L}_1\mathscr{L}_1\mathscr{L}_1\mathscr{L}_2-\mathscr{L}_2\mathscr{L}_1\mathscr{L}_1\mathscr{L}_1+3\mathscr{L}_1\mathscr{L}_2\mathscr{L}_1\mathscr{L}_1-3\mathscr{L}_1\mathscr{L}_1\mathscr{L}_2\mathscr{L}_1)\Delta  t^4\\
&+(c_1d_1-\frac{3}{8}d_1-\half{c_1d_1^2}-\half{c_1^2d_1}+\frac{1}{4}d_1^2+\frac{1}{16})\\
&\quad \times (\mathscr{L}_2\mathscr{L}_2\mathscr{L}_1\mathscr{L}_1-\mathscr{L}_1\mathscr{L}_1\mathscr{L}_2\mathscr{L}_2+2\mathscr{L}_1\mathscr{L}_2\mathscr{L}_1\mathscr{L}_2-2\mathscr{L}_2\mathscr{L}_1\mathscr{L}_2\mathscr{L}_1)\Delta  t^4\\
&+(\half{c_1d_1}-\frac{1}{8}{d_1}-\half{c_1^2d_1}+\frac{1}{48})\\
&\quad \times (\mathscr{L}_1\mathscr{L}_2\mathscr{L}_2\mathscr{L}_2-\mathscr{L}_2\mathscr{L}_2\mathscr{L}_2\mathscr{L}_1+3\mathscr{L}_2\mathscr{L}_2\mathscr{L}_1\mathscr{L}_2-3\mathscr{L}_2\mathscr{L}_1\mathscr{L}_2\mathscr{L}_2)\Delta  t^4+\cdots
 \end{aligned}
\end{equation}

Let us define,
 \begin{equation}\label{eq: D4}
\begin{aligned}
\mathscr D_{41}\triangleq& \mathscr{L}_1\mathscr{L}_1\mathscr{L}_1\mathscr{L}_2-\mathscr{L}_2\mathscr{L}_1\mathscr{L}_1\mathscr{L}_1+3\mathscr{L}_1\mathscr{L}_2\mathscr{L}_1\mathscr{L}_1-3\mathscr{L}_1\mathscr{L}_1\mathscr{L}_2\mathscr{L}_1,\\
\mathscr D_{42}\triangleq& \mathscr{L}_2\mathscr{L}_2\mathscr{L}_1\mathscr{L}_1-\mathscr{L}_1\mathscr{L}_1\mathscr{L}_2\mathscr{L}_2+2\mathscr{L}_1\mathscr{L}_2\mathscr{L}_1\mathscr{L}_2-2\mathscr{L}_2\mathscr{L}_1\mathscr{L}_2\mathscr{L}_1,\\
\mathscr D_{43}\triangleq& \mathscr{L}_1\mathscr{L}_2\mathscr{L}_2\mathscr{L}_2-\mathscr{L}_2\mathscr{L}_2\mathscr{L}_2\mathscr{L}_1+3\mathscr{L}_2\mathscr{L}_2\mathscr{L}_1\mathscr{L}_2-3\mathscr{L}_2\mathscr{L}_1\mathscr{L}_2\mathscr{L}_2.
\end{aligned} 
\end{equation}

A similar analysis yields,
\begin{theorem}
The operators in \eqref{eq: D4} are of the orders,
\begin{equation}\label{eq:thm1}
\begin{aligned}
&1)\mathscr \sD_{41}=\mathcal{O}(\frac{1}{\veps^4}),\\
&2)\mathscr \sD_{42}=\mathcal{O}(\frac{1}{\veps^2}),\\
&3)\mathscr \sD_{43}=\mathcal{O}(1).
\end{aligned}
\end{equation}
\end{theorem}

%Based on these estimates, we have for the original impulse method:
%\begin{coro}
% When $c_1=0, c_2=c_3=\half, c_4=0, d_1=0, d_2=1, d_3=d_4=0$, we have 
% \begin{equation}
% \widetilde{\sR}(c_1,c_2,c_3,c_4;d_1,d_2,d_3,d_4)=1+  \cal{O}(\Delta  t^2)  + \cal{O}(\frac{\Delta  t^3}{\veps^2}) + \mathcal{O}(\frac{\Delta  t^4}{\veps^4})
% + \cdots 
%\end{equation}
%\end{coro} 
\smallskip
In order to minimize the error, especially the first few terms of $ \widetilde{\sR}$, we will choose the parameters so that coefficients of $\mathscr D_{41}$ and $\mathscr D_{42}$ are zero  (Notice that the coefficient of $\mathscr D_{31}$ is twice of the one of $\mathscr D_{41}$, so it will be automatically zero). This leaves us with the nonlinear equations for $c_1$ and $d_1$:
\begin{equation}
 \begin{aligned}
  \half{c_1d_1}-\frac{1}{4}d_1-\half{c_1d_1^2}+\frac{1}{4}d_1^2+\frac{1}{24}=0,\\
  c_1d_1-\frac{3}{8}d_1-\half{c_1d_1^2}-\half{c_1^2d_1}+\frac{1}{4}d_1^2+\frac{1}{16}=0.
 \end{aligned}
\end{equation}
These equations can be simplified. In particular, $d_1$ is the only root of $6z^3-12z^2+6z-1,$ and $c_1=\half d_1.$

\begin{coro}
 When $c_1=\frac{\sqrt[3]{2}}{6}+\frac{\sqrt[3]{4}}{12}+\frac{1}{3}, d_1=2c_1$, the coefficients of $\sD_{31},$ $\sD_{32},$ $\sD_{41},$ $\sD_{42},$ and $\sD_{43}$ in \eqref{eq: dt34} are all zeros.
\end{coro}

The elimination of the $\sD_{32}$ and $\sD_{43}$ terms seem to be a coincidence.  This method will be referred to as {\bf impulse IV.} With a direct calculation, one can show that
$c_1=\frac{1}{2\big(2-2^{\frac13}\big)},$ which surprisingly,  coincides with the coefficients of the well known 4th order symplectic integrator \cite{yoshida1990construction}. Of course, the symplectic method in \cite{yoshida1990construction} is based on the splitting of the kinetic and potential energy, while our splitting is between the fast and slow forces. 

\subsection{A Numerical Test:  A nonlinear coupled oscillator}
Before we look further into the properties of the splitting methods, we present some numerical results for
a nonlinear oscillator problem \cite{leimkuhler2004simulating}, governed by the equations,
\begin{equation}
 \left\{
  \begin{aligned}
  \ddot{q}= &\frac{1-\|q\|}{\|q\|} q + \beta \|\theta-q\|^2 (\theta -q)  + \frac{1}{\veps^2}  (\theta -q), \\
 \ddot{\theta}= & -\beta \|\theta-q\|^2 (\theta -q)  - \frac{1}{\veps^2}  (\theta -q). 
 \end{aligned}
 \right.
\end{equation}
Here $q, \theta \in \mathbb{R}^2$. This is a Hamiltonian system with  potential energy given by,
\begin{equation}
V(q,\theta)= \frac{1}{2\veps^2}\|\theta - q\|^2 + \frac{\beta}{4}\|\theta-q\|^4 + \half\big( \|q\| -1\big)^2.
\end{equation}
In our tests, the parameters are chosen as follows: $\veps=0.1,$ $\beta=0.1$, $q(0)=(1,0),$ $\theta(0)=(1.01,0),$
$\dot{q}(0)=(0,1),$ $\dot{\theta}(0)=(0,0.05)$, $\dt=0.01$, and $\Dt=0.12.$  

\begin{figure}[htp]
\begin{center}
\includegraphics[scale=0.55]{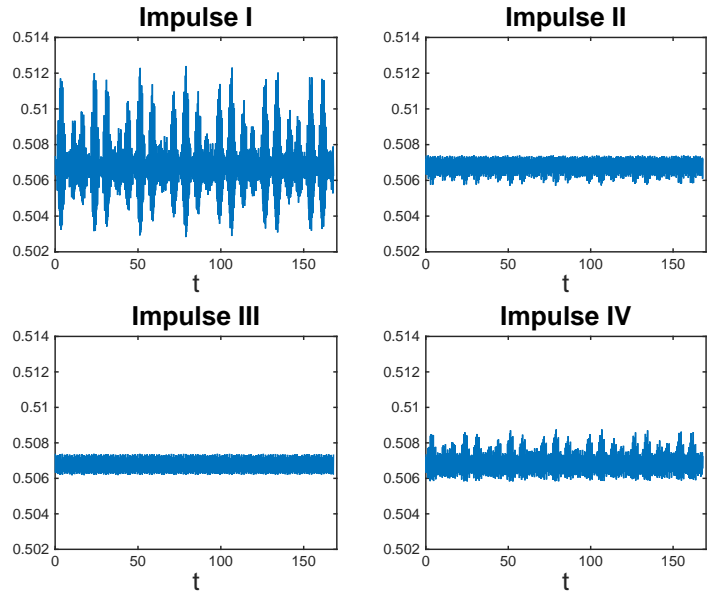}
\caption{Comparison of the energy conservation for the four methods.}
\label{fig: H}
\end{center}
\end{figure}

First, we show the total energy computed from each method in Fig. \ref{fig: H}.  We find that the new impulse methods have much better performance in the energy conservation: The fluctuation is much smaller than the original impulse method (impulse {\bf I}). For problems where the energy is more relevant than the actually trajectories, e.g., producing various statistical ensembles, the new methods seem to be more promising. Among the new impulse methods, the method \biii { seems} to have the best results.   It is clear, however, much deeper analysis is needed to understand the accuracy of the methods toward computing different quantities. 
We will present some preliminary analysis in the next section.

\begin{figure}[htp]
\begin{center}
\includegraphics[scale=0.55]{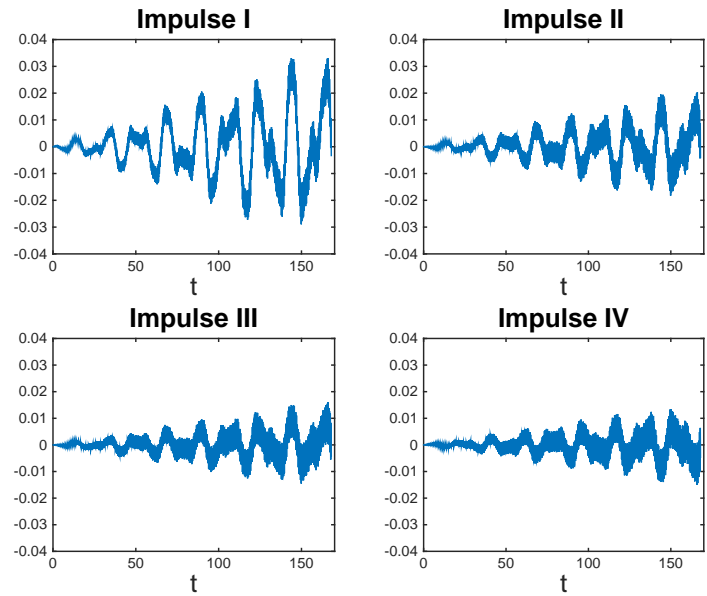}
\caption{Comparison of the error in $q_1$ for the four methods.}
\label{fig: q1}
\end{center}
\end{figure}

In Fig. \ref{fig: q1}, we show the error of $q_1$ for the numerical approximations obtained from the impulse methods. 
The error is estimated by comparing the approximate solutions to a solution computed with very small step size. We observe that the accuracy is gradually improved for the impulse methods {\bf I} to {\bf IV}. However, the error would grow in all cases, which can be attributed to the Lyapunov instability inherent in most Hamiltonian systems.

\begin{figure}[htp]
\begin{center}
\includegraphics[scale=0.55]{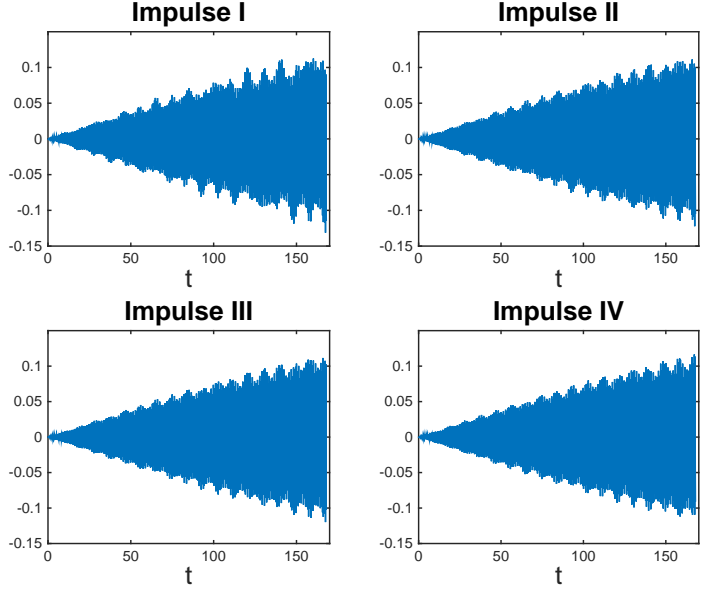}
\caption{Comparison of the error in $p_1$ for the four methods.}
\label{fig: p1}
\end{center}
\end{figure}

In Fig. \ref{fig: p1}, we show the momentum computed from the four methods. Surprisingly, they exhibit similar accuracy, and the new impulse methods show little improvement. 

\subsection{Some preliminary analysis}

\subsubsection{Energy conservation of multiscale Hamiltonian systems}
For those ODEs \eqref{eq:hamilton} that come from Hamiltonian systems, i.e.,  $f=-M^{-1}\nabla W$ and $g=-M^{-1}\nabla V,$
 we define two Hamiltonians that correspond to the splitting of the operator $\sL$, 
\begin{equation}
\begin{aligned}
 H_2= &W(x),\\
 H_1= &\frac{1}{\veps^2} V(x) + \half p^T M^{-1} p.
\end{aligned}
\end{equation}
Here $M$ is the mass matrix, and $p=Mv$ is the momentum variable.

Due to the approximation, the energy associated with the dynamical system will not be exactly conserved. However, one of the celebrated results in geometric integrators is that an approximate Hamiltonian often exists, and it is conserved exactly by the numerical method \cite{feng1986difference,hairer2006geometric,leimkuhler2004simulating,ruth1983canonical,yoshida1993recent}. 

For impulse \bi the energy conservation property can be analyzed using the backward analysis  \cite{hairer2006geometric,yoshida1993recent}, which asserts that the approximate solution is a more accurate solution of another Hamiltonian system with a Hamiltonian $H_S$, known as the {\it shadow Hamiltonian}. For the impulse method {\bf I}, the analysis shows that the shadow Hamiltonian, up to the order $\Delta t^2$, is given by \cite{hairer2006geometric},
\begin{equation}
 H_S = H_1+H_2 + \frac{\Delta t^2}{12}\{\{H_2,H_1\},H_1\} - \frac{\Delta t^2}{24}\{H_2,\{H_2,H_1\}\}.
\end{equation}
This is because the operator approximation can be written as,
\begin{equation}
 e^{\half\Delta  t \mathscr{L}_2}e^{\Delta  t \mathscr{L}_1} e^{\half\Delta  t\mathscr{L}_2}
 = e^{\sL \Dt + \frac{\Dt^3}{12} \sD_{31} -  \frac{\Dt^3}{24} \sD_{32} + \cdots}.
\end{equation}
Here, $\sD_{31}=[[\sL_2,\sL_1],\sL_1]$ and $\sD_{32}=[\sL_2,[\sL_2,\sL_1]]$; $[\quad]$ and $\{\quad\}$ stand for the commutator (Lie derivative) and Poisson bracket, respectively. 
 
 In particular, we have that,
\begin{equation}
 \begin{aligned}
 \{H_2,H_1\}=&-p^T M^{-1 }\nabla W(x),\\
  \{\{H_1,H_2\},H_2\} = & p^TM^{-1 }\nabla^2WM^{-1 }p-\frac{1}{\veps^2}\nabla W^TM^{-1 }\nabla V(x),  \\
   \{\{H_2,H_1\},H_1\}=&\nabla W^T M^{-1 } \nabla W.  
\end{aligned}
\end{equation}

As a result, the conservation of the energy at this level is dominated by the $\mathcal{O}\big(\frac{\Delta t^2}{\veps^2}\big)$ term, which is due to the presence of $\sD_{31}$ in the error $\wt{\sR}$. 

In contrast, the same calculation for the impulse method \bii yields,
\begin{equation}
 H_S = H  + \frac{17\Delta t^2}{96}\{H_2,\{H_2,H_1\}\}=H+\mathcal{O}\big({\Delta t^2}\big)+\cdots
\end{equation}
As a result, the better energy conservation can be attributed to the elimination of the $\sD_{31}$ term in the error.

\subsubsection{Resonance instability}

\begin{figure}[htbp]
\begin{center}
\includegraphics[scale=0.5]{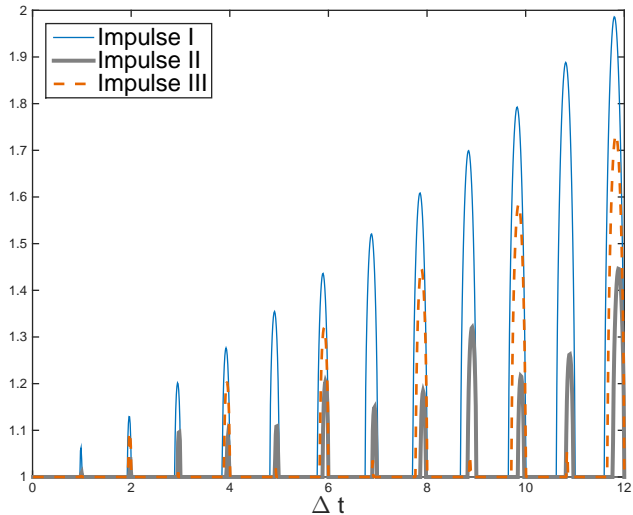}
\caption{Spectral radius of the propagation matrix.}
\label{fig: eig}
\end{center}
\end{figure}

Another outstanding issue raised by previous works is the resonance, which occurs for certain choices of the slow time step $\Dt$  \cite{schlick2010molecular,barth1998extrapolation}. Following the analysis in \cite{barth1998extrapolation,barth1998overcoming}, and in  particular the example in \cite{schlick2010molecular}, we consider a scalar problem where 
$f(x)= -\big(\frac{\pi}{5}\big)^2 x$ and $g(x)= -\pi^2 x.$ In Fig. \ref{fig: eig}, we show the spectral radius of the propagation matrix. We observe that the new impulse methods exhibit similar resonance phenomena: When the large step size $\Dt$ is around an integer multiple of half of the period (T=2) associated with the fast scale, instability occurs. 

\section{Another example: Dynamics of Octane}

Here, we consider the Octane molecule with 26 atoms.  The impulse methods have been implemented within TINKER \cite{ponder1987tinker}. To properly quantify the error, the `exact' solution is represented by the solution computed with the Verlet's method with small step size $10^{-5}ps.$ In the impulse method, we choose $\Dt=2.4\times 10^{-4}ps$ and $\dt=\Dt/24.$ All the simulations are conducted for $6ps$ period.

\begin{figure}[htbp]
\begin{center}
\includegraphics[scale=0.5]{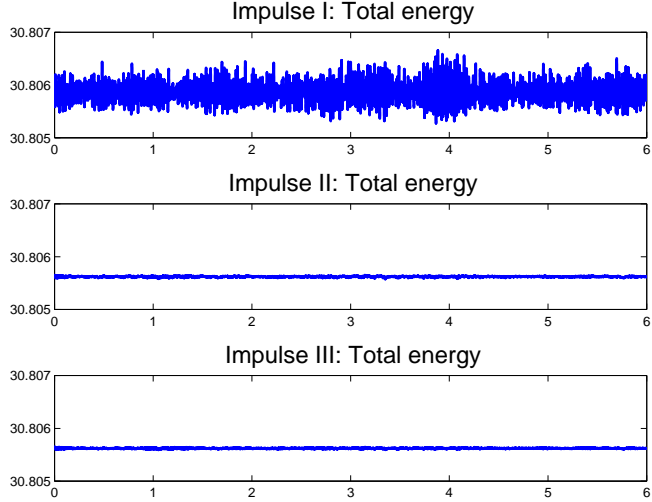}
\caption{Comparison of the energy conservation.}
\label{fig: oct-eng}
\end{center}
\end{figure}

In Fig.   \ref{fig: oct-eng}, we show the total energy computed from impulse methods \bi to \biii. Again we observe that the new impulse methods have much less fluctuation of the energy, indicating a better energy conservation property.

\begin{figure}[htbp]
\begin{center}
\includegraphics[scale=0.5]{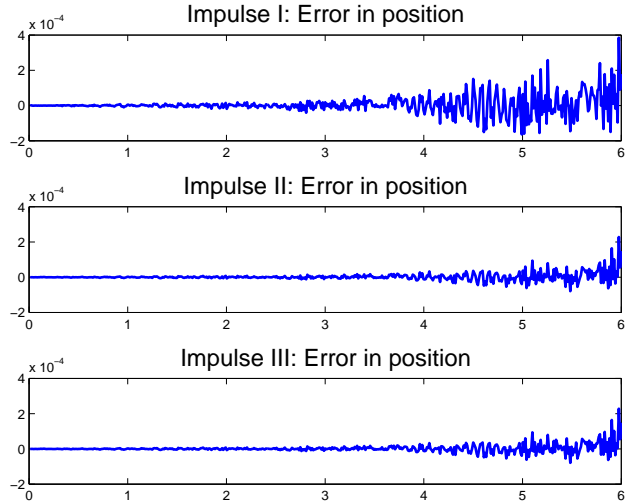}
\caption{Comparison of the error in $x_1$.}
\label{fig: oct-pos}
\end{center}
\end{figure}

\begin{figure}[htbp]
\begin{center}
\includegraphics[scale=0.5]{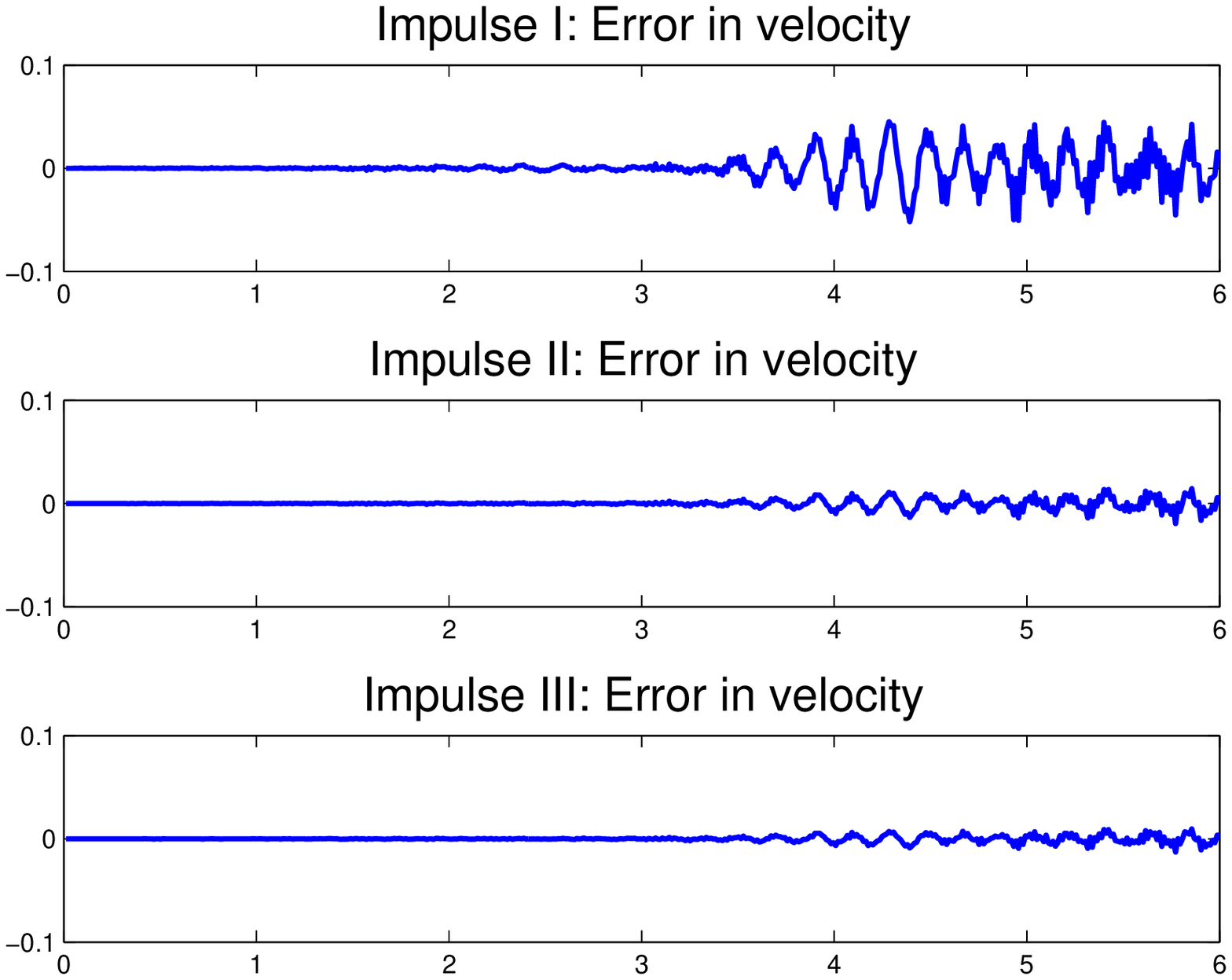}
\caption{Comparison of the error in $v_1$.}
\label{fig: oct-vel}
\end{center}
\end{figure}

Next we look at the error in the position of the first atom (first component). The results are shown in Fig. \ref{fig: oct-pos}. There are some improvement of the accuracy from impulse methods \bi to impulse \biii.
Such improvement is also observed in the velocity (first component $v_1$), as can be seen in Fig. \ref{fig: oct-vel}.

\section{Summary and discussion}

We have developed some new impulse methods for the numerical approximation of molecular systems with multiple time scales, which are represented by interactions of different magnitude. Motivated by the operator-splitting approach of the original impulse method \cite{grubmuller1991generalized,tuckerman1992reversible}, we sought general splitting methods that involve more fractional steps. A novel aspect in our approach is the systematic procedure for finding an expansion of the error, which in turn sheds light on the selection of the coefficients so that the accuracy can be improved.

For multiscale ODEs, our analysis revealed that the terms in the error can depend on both the large time step $\Dt$,  and $\veps$, which represents the separation of the scales. Based on the order of the first few terms, we choose the coefficients so that the terms with the largest magnitude are eliminated. This leads to several splitting methods that can be viewed as generalized impulse methods. Numerical tests have been
conducted, which have confirmed the improved accuracy. The biggest improvement has been observed in the energy conservation. This has been analyzed with some preliminary study of the modified Hamiltonian. It is thus expected that the new methods would produce better results when the micro-canonical ensemble distribution is of interest. 

This approach, however, is by no means complete. From a practical viewpoint, the Hamiltonian system \eqref{eq:hamilton} models a system in isolation. In practice, often of interest are extended systems, where the external conditions are modeled by introducing additional variables, e.g., heat and pressure bath, or by introducing stochastic forces, e.g., the Langevin dynamics. Operator-splitting methods have been widely used for  extended systems, e.g., in \cite{martyna1996explicit}, and the impulse methods have also been applied to Langevin dynamics \cite{skeel2002impulse} as well. Extending the current framework to those problems might produce new integrators with other capabilities, and it will be explored in our future works. 

The problem considered in this work belongs to stiff ODEs, for which many numerical methods have been developed, e.g., implicit 
Runge-Kutta methods and BDF methods, and they can be found in standard textbooks  \cite{deuflhard2002scientific}. Meanwhile,
there have been significant recent  progress in developing efficient computational methods for  such 
 dynamical systems with multiple time scales, e.g., the heterogeneous multiscale method (HMM) \cite{engquist2005heterogeneous,ariel2009multiscale,brumm2013heterogeneous,fatkullin2004computational}, the equation-free method \cite{kevrekidis2003equation}, the FLAVOR method \cite{tao2010nonintrusive}, the reversible averaging integrator \cite{leimkuhler2004simulating}, etc. These methods demonstrate resemblance to the impulse method in that they introduce multiple time steps ($\dt$ and $\Dt$) to capture the multiple scales. On the other hand, an averaging procedure is usually involved on the fastest time scale to compute an effective force on quantities that evolve on the slow time scale.  In addition, some of these methods 
 assume the existence and explicit form of slow variables.   At this point, we are not aware of the application of these methods to biomolecular modes.

\section*{Acknowledgement} This project was completed when Yuan was participating in the MASS program in the Department of Mathematics at Penn State University in the fall of 2014. She would like to acknowledge the support from her home institution, Wuhan University, and  the MASS program for the research opportunity.  She would also like to thank Xiaojie Wu for the help with the computing facility at Penn State. One of the test problems was from the software TINKER \cite{ponder1987tinker}. 
\bibliographystyle{elsarticle-num}
\bibliography{msode}

\end{document}